\documentclass[reqno]{amsart}

\usepackage{mathtools}
\usepackage{amsmath}

\usepackage[dvipsnames]{xcolor}
\usepackage
[colorlinks=true,linkcolor=Maroon,citecolor=OliveGreen,backref]
{hyperref}

\newtheorem{theorem}{Theorem}[section]

\theoremstyle{definition}
\newtheorem{remark}[theorem]{Remark}

\usepackage{amssymb}
\newcommand{\eps}{\epsilon}

\def\P{\mathbf{P}}
\def\E{\mathbf{E}}

\def\calC{\mathcal{C}}

\newcommand\br[1]{\left({#1}\right)}
\newcommand\floor[1]{\left\lfloor{#1}\right\rfloor}

\newcommand{\refmain}[1]{\cite[#1]{EG1}}

\begin{document}

\title{Probability of generation by random permutations of given cycle type}

\author{Sean Eberhard}
\address{Centre for Mathematical Sciences\\
    Wilberforce Road\\
    Cambridge CB3~0WB (U.K.)}
\email{eberhard@maths.cam.ac.uk}

\author{Daniele Garzoni}
\address{School of Mathematical Sciences, Tel Aviv University\\ Tel Aviv 69978 (Israel)}
\email{danieleg@mail.tau.ac.il}


\begin{abstract}
Suppose $\pi$ and $\pi'$ are two random elements of $S_n$ with constrained cycle types
such that $\pi$ has $x n^{1/2}$ fixed points and $yn/2$ two-cycles,
and likewise $\pi'$ has $x' n^{1/2}$ fixed points and $y'n/2$ two-cycles.
We show that the events that $G = \langle \pi, \pi' \rangle$ is transitive and $G \geq A_n$
both have probability approximately
\[
    (1 - yy')^{1/2} \exp\br{- \frac{xx' + \frac12 x^2 y' + \frac12 {x'}^2 y}{1 - yy'}},
\]
provided $(x, x')$ is not close to $(0, \infty)$ or $(\infty, 0)$.
This formula is derived from some preliminary results in a recent paper of the authors.
As an application, we show that two uniformly random elements of uniformly random conjugacy classes of $S_n$ generate the group with probability about 51\%.
\end{abstract}

\maketitle

\section{Introduction}

Consider a sequence of pairs of conjugacy classes $\calC_n, \calC'_n$ of $S_n$ (one pair for each positive integer $n$).
Let $c_i$ and $c'_i$ be the number of $i$-cycles of an element of $\calC_n$ and $\calC'_n$, respectively.
We are interested in the group $G$ generated by uniformly random $\pi \in \calC_n$ and $\pi' \in \calC'_n$.

\begin{theorem}
    \label{thm:main}
    Assume
    \begin{align*}
        c_1 / n^{1/2} &\to x, & 2 c_2 / n &\to y,\\
        c'_1 / n^{1/2} &\to x', & 2 c'_2 / n &\to y',
    \end{align*}
    where $0 \leq x, x' \leq \infty$ and $0 \leq y, y' \leq 1$.
    If $yy' < 1$, assume $(x, x') \notin \{(0, \infty), (\infty, 0)\}$.
    Then
    \begin{equation}
        \label{eq:main-estimate}
        \P(G \geq A_n) \to (1 - yy')^{1/2} \exp\br{-\frac{xx' + \frac12 x^2 y' + \frac12 {x'}^2 y}{1 - yy'}}.
    \end{equation}
    The right-hand side is interpretted (continuously) as zero if $y = y' = 1$ or $x = \infty$ or $x' = \infty$.
    The same estimate holds for $\P(G~\text{transitive})$.
\end{theorem}

This theorem essentially generalizes \refmain{Theorem~1.1},
which describes the boundary cases of \eqref{eq:main-estimate}.
We will deduce Theorem~\ref{thm:main} from some of the preliminary results of \cite{EG1}.

\begin{remark}
    In the indeterminate case $(x, x') \in \{(0, \infty), (\infty, 0)\}$, $\P(G \geq A_n)$ can be close to 0 or to 1.
    Some examples are given in \refmain{Section~2.5}.
\end{remark}

\section{Proof}

By \refmain{Theorem~1.1(2)} we may assume $x + x' < \infty$ and $yy' < 1$.
By \refmain{Lemma~4.5 and Lemma~4.9}, $\P(G~\text{transitive})$ and $\P(G \ge A_n)$ differ by $e^{-\Omega(n)}$, so we may focus on $\P(G~\text{transitive})$.
As in \refmain{Section~3} let $N_k$ be the number of orbits of $G$ of size $k$ and let
\[
    N = \sum_{k=1}^{\floor{n/2}} N_k.
\]
Then $G$ is transitive if and only if $N = 0$.
By \refmain{Lemma~3.3 and Lemma~3.4}, $\E N_k \leq \eps_k$ for some constant $\eps_k$ such that $\sum_k \eps_k < \infty$, which allows us to apply the dominated convergence theorem.
Now fix $k$.
By \refmain{(3)},
\[
    \E N_k = \sum \binom{n}{k}^{-1} \binom{c_1}{d_1} \cdots \binom{c_k}{d_k} \binom{c'_1}{d'_1} \cdots \binom{c'_k}{d'_k} p(d_1, \dots, d_k ; d'_1, \dots, d'_k),
\]
where the sum goes over all solutions to
\begin{align*}
    k &= d_1 1 + d_2 2 + \cdots + d_k k \\
    k &= d'_1 1 + d'_2 2 + \cdots + d'_k k
\end{align*}
such that $0 \leq d_i \leq c_i$ and $0 \leq d'_i \leq c'_i$ for each $i$,
and where $p(d_1, \dots, d_k ; d'_1, \dots, d'_k)$ is the probability that
random permutations $\tau, \tau' \in S_k$ with cycle types $1^{d_1} \cdots k^{d_k}$ and $1^{d'_1} \cdots k^{d'_k}$
generate a transitive subgroup.
Write
\[
    \E N_k = \Sigma_1 + \Sigma_2,
\]
where $\Sigma_1$ is the sum of all terms with $d_i = d'_i = 0$ for each $i \geq 3$
and $\Sigma_2$ is the sum of all other terms.
If $d_i + d'_i > 0$ for any $i \geq 3$ then
\[
    \binom{c_1}{d_1} \cdots \binom{c_k}{d_k} \binom{c'_1}{d'_1} \cdots \binom{c'_k}{d'_k}
    \leq O_k(n^{d_1/2 + d_2 + \cdots + d_k} n^{d'_1/2 + d'_2 + \cdots + d'_k})
    \leq O_k(n^{k - 1/2}).
\]
Hence
\[
    \Sigma_2 \leq O_k(n^{-1/2}).
\]
Now consider $\Sigma_1$.
We have $p(d_1, d_2 ; d'_1, d'_2) = 0$ unless $d_1 + d'_1 \leq 2$, so if $k = 2m + 1$ is odd we have
\[
    \Sigma_1 = \binom{n}{k}^{-1} \binom{c_1}{1} \binom{c_2}{m} \binom{c'_1}{1} \binom{c'_2}{m} p(1, m ; 1, m),
\]
while if $k = 2m$ is even then $\Sigma_1$ is the sum of three terms:
\begin{align*}
    \Sigma_1
    &=\binom{n}{k}^{-1} \binom{c_1}{2} \binom{c_2}{m - 1} \binom{c'_2}{m} p(2, m - 1; 0, m)\\
    &+ \binom{n}{k}^{-1} \binom{c_2}{m} \binom{c'_1}{2} \binom{c'_2}{m-1} p(0, m ; 2, m - 1) \\
    &\qquad+ \binom{n}{k}^{-1} \binom{c_2}{m} \binom{c'_2}{m} p(0, m ; 0, m).
\end{align*}
The three calculations
\begin{align*}
    p(1, m; 1, m) &= \frac{(2m+1)!}{((2m+1)! / 2^m m!)^2},\\
    p(2, m-1; 0, m) &= \frac{(2m)! / 2}{((2m)! / 2^{m-1} 2! (m-1)!) \times ((2m)! / 2^m m!)},\\
    p(0, m; 0, m) &= \frac{(2m-1)!}{((2m)! / 2^m m!)^2}
\end{align*}
can be left to the reader (cf.~\refmain{Lemma~2.1}).
Hence, as $n \to \infty$ with $k$ fixed,
\[
    \Sigma_1 \to
    \begin{cases}
    x y^m x' {y'}^m &: k = 2m+1 > 0, \\
    (x^2 / y + {x'}^2 / y') y^m {y'}^m / 2
    + y^m {y'}^m / (2m) &: k = 2m > 0.
    \end{cases}
\]
It follows by the dominated convergence theorem that
\begin{align*}
    \E N
    &= \sum_{k=1}^{\floor{n/2}} \E N_k \\
    &\to \sum_{k = 2m+1 > 0} x x' (yy')^m + \sum_{k = 2m > 0} \br{(x^2 / y + {x'}^2 / y') (yy')^m / 2 + (yy')^m / (2m)} \\
    &= (xx' + x^2 y'/2 + {x'}^2 y/2) / (1 - yy') - \frac12 \log(1 - yy').
\end{align*}
Hence,
by \refmain{Theorem~3.10},
\[
    \P(N = 0) = e^{-\E N} + o(1) \to (1 - yy')^{1/2} \exp\br{- \frac{xx' + \frac12 x^2 y' + \frac12 {x'}^2 y}{1- yy'}}.
\]
This completes the proof.

\section{Application}

Let $\calC_n$ and $\calC'_n$ be uniformly random conjugacy classes of $S_n$,
and let $\pi \in \calC_n$ and $\pi' \in \calC'_n$ be uniformly random.
Let $G = \langle \pi, \pi'\rangle$.
What is the probability that $G \geq A_n$?
The Hardy--Ramanujan asymptotic for the partition function states that
\[
    p(n) \sim \frac{a}{n} \exp(2 b n^{1/2}),
\]
where $a = 1 / (4\sqrt{3})$ and $b = \pi/\sqrt{6}$ (see \cite[Proposition~VIII.6, p.~578]{flajolet-sedgewick}).
Since the cycle type of $\pi$ is that of a uniformly random partition of $n$, it follows that
\begin{align*}
    \P(c_1 \geq x n^{1/2}, c_2 \geq y n^{1/2})
    &\sim \frac{n}{n - (x + 2y)n^{1/2}} \exp(2b(n - (x + 2y)n^{1/2})^{1/2} - 2b n^{1/2}) \\
    &\to \exp(- b (x+2y)).
\end{align*}
In particular note that $c_2 = o(n)$ with high probability.
Hence Theorem~\ref{thm:main} implies that $\P(G \ge A_n)$ and $\P(G~\text{transitive})$ are both asymptotically
\[
    b^2 \int_0^\infty \int_0^\infty \exp(-xx' - b(x + x')) \, dx \, dx' = b^2 e^{b^2} \int_{b^2}^\infty t^{-1} e^{-t} \, dt \approx 0.6889.
\]
A random conjugacy class consists of even permutations with probability tending to $1/2$.
Hence $\P(G = A_n) \approx 0.1722$ and $\P(G = S_n) \approx 0.5167$.

\bibliographystyle{amsplain-ac}
\bibliography{refs}
\end{document}